\documentclass{amsart}
\usepackage{amssymb, amscd}
\usepackage[all]{xy}
\usepackage{epic}
\usepackage{url}
\usepackage[english]{babel}

\makeindex

\makeatletter
\def\@doGroups@#1#2#3{#1{#3}#2{#3}}
\def\doGroups#1#2#3#4{%
\def\@doGroups{\@ifnextchar#1\@doGroups@@{\@doGroups@#3#4}}%
\def\@doGroups@@#1##1#2{#3{##1}%
\@ifnextchar#1\@doGroups@@{#4}}%
\@doGroups}

\let\@DefHandle\index

\newcommand\Definition{\doGroups[]\@DefHandle\emph}

\makeatother

\numberwithin{equation}{section}

\newtheorem{theorem}{Theorem}[section]
\newtheorem{lemma}[theorem]{Lemma}
\newtheorem{proposition}[theorem]{Proposition}
\newtheorem{corollary}[theorem]{Corollary}

\newtheorem{definition-lemma}[theorem]{Definition-Lemma}

\theoremstyle{definition}

\theoremstyle{remark}

\newtheorem{property}[theorem]{Property}

\def\AA{{\mathbb A}}

\def\FF{{\mathbb F}}

\def\PP{{\mathbb P}}

\def\cO{{\mathcal O}}

\makeatletter
\newcommand\symb[1]{\mathop {\operator@font #1}\nolimits}
\makeatother

\newcommand\riso{\mathrel{\hskip2pt\raise-2.5pt\hbox{$\widetilde{\phantom{xx}}$}\kern-16pt\longrightarrow}}

\newcommand\disjunion{\setbox0\hbox{$\cup$}%
\dimen0\wd0
\setbox1=\hbox to\dimen0{\hfil.\hfil}
\dimen0=\dp1
\mathbin{\vbox{\offinterlineskip\box1\box0}}}

\newcommand\Disjunion{\setbox0\hbox{$\bigcup$}%
\dimen0\wd0
\setbox1=\hbox to\dimen0{\hfil.\hfil}
\dimen0=\dp1
\mathbin{\vbox{\offinterlineskip\box1\box0}}}

\newcommand\proofsquare{\nobreak\hfill \hbox{%
\vrule height 5pt
\kern-.4pt
 \vbox{%
\hrule width 5pt depth0pt height.4pt
 \kern4.6pt \hrule  }
\kern-3.75pt
\vrule height 5pt}\kern1pt
\par}

\countdef\Yxpos255
\countdef\Yypos254
\countdef\Yxht 253
\countdef\Yxold252
\def\oneYoung#1{\ifnum#1=0
\put(\Yxpos,\Yypos){\line(0,1){0.2}\line(1,0){1}}%
\else
\Yxht=#1
{\ifnum#1<\Yxold
\Yxht=\Yxold
\fi
\put(\Yxpos,\Yypos){\line(0,1){\Yxht}}}%
\advance\Yxht by 1
\multiput(\Yxpos,\Yypos)(0,1){\Yxht}{\line(1,0){1}}%
\fi
\Yxold=#1
}

\def\Young(#1,#2)(#3){%
\def\nxt##1,{\advance\Yxpos by 1
\ifx!##1\else\oneYoung{##1}\expandafter\nxt\fi}%
\Yxold=0
\Yxpos=0
\advance\Yxpos by -1
\Yypos=0
\put(#1,#2){\nxt#3,!,
\put(\Yxpos,\Yypos){\line(0,1){\Yxold}}}}

\begin{document}

\title[The Galois Closure of the Garcia-Stichtenoth Tower]{
The Galois Closure of the Garcia-Stichtenoth Tower}

\author{Alexey Zaytsev}

\address{Faculteit Wiskunde en Informatica, University of
Amsterdam, Plantage Muidergracht 24, 1018 TV Amsterdam, The Netherlands}

\email{azaitsev@science.uva.nl}


\begin{abstract}
We describe the Galois closure of the Garcia-Stichtenoth tower and prove that
it is optimal.
\end{abstract}
\maketitle

\begin{section}{Introduction}
In 1996 Garcia and Stichtenoth constructed in \cite{G-S1} a tower
of Artin-Schreier covers
$$
\ldots \rightarrow X_{i} \rightarrow X_{i-1} \rightarrow \ldots
X_{1} \rightarrow \PP^{1}
$$
which are defined over the finite field $\FF_{q^2}$ and given by a simple recursive
equation such that
$$
\lim_{n \rightarrow \infty}N(X_n)/g(X_n)=q-1,
$$
where $N(X_n)$ is the number of $\FF_{q^2}$-rational points and $g(X_n)$ is the genus of $X_n$.
In this note we construct the Galois closure of this tower, i.e.
the tower of covers
$$
\ldots \rightarrow \tilde{X}_{i} \rightarrow \tilde{X}_{i-1}
\rightarrow \ldots \tilde{X}_{1} \rightarrow \PP^{1}
$$
such that $\tilde{X}_i$ is the Galois closure of the cover $X_i
\rightarrow \PP^1$. We give explicit formulas of the genus, estimate
the number of $\FF_{q^2}$-rational points of the curve $\tilde{X}_i$ and show that the tower is optimal as well,
i.e., it reaches the Drinfeld-Vl\u{a}du\c{t} upper-bound.
\end{section}

\begin{section}{Generators for the Galois Closure}
Let $p$ be an odd prime number and $K=\FF_{p^2}$ be a finite field of cardinality $p^2$.
Garcia and Stichtenoth described in \cite{G-S1} a tower of curves or function
fields over $K$ by defining recursively  fields $T_m:=K(x_1,\ldots,x_m)$
with $x_{i+1}$ satisfying the equation
$$
x_{i+1}^p+x_{i+1}=x_i^{p+1}/(x_{i}^{p}+x_{i}). \eqno(1)
$$
We shall write $\wp(x)$ for the expression $x^p+x$ and we let $g$ be the
rational function $x^{p+1}/(x^p+x)$ in $\FF_p(x)$. Furthermore, we set
$h=(x^{p-1}-1)/(x^{p-1}+1) \in \FF_p(x)$.
Then $T_n$ is obtained
from $T_{n-1}$ by adjoining a root $y=x_n$ of the equation $\wp(y)=g(x_{n-1})$.
We let $\tilde{T}_n$ be the Galois closure of $T_n$ over $T_1$ and
$\Gamma_n$ the Galois group of $\tilde{T}_n$ over $T_1$.

We set $K_{-}=\{ \alpha \in K : \alpha^p=-\alpha \}$. If $c=(c_1,\ldots,c_n)
\in K_{-}^n$ (for $n \geq 2$) then we denote by $u_c$ a root of
$$
f_c:=X^p+X-g(u_{c^{\prime}}+c_n),
$$
where $c^{\prime}$ is the shortened vector $c^{\prime}=(c_1,\ldots,c_{n-1})$ and
for $n=1$, the element $u_{c_1}$ is a root of polynomial
$$
f_{c_1}:=X^p+X-g(x_2+c_1).
$$
We make the generators more precise in the following proposition.
\begin{proposition} For $n \geq 3$
the field $\tilde{T}_n$ is generated over $\tilde{T}_{n-1}$ by
adjoining all the elements $u_c$ with $c \in K_{-}^{n-2}$.
\begin{proof}
The proof is by induction on $n$. If $n=3$ then the field $\tilde{T}_3$ is
the composite of the fields $\tilde{T}_2(\sigma(x_3))$ with $\sigma$
 running through $\Gamma_3$. By applying $\sigma$ to (1) one sees that
$\sigma(x_2)=x_2+c_1$ for some $c_1 \in K_{-}$. Similarly, one observes that
$\wp(x_3)=g(x_2+c_1)$, i.e.,
$\sigma(x_3)=u_{c_1}+c_2$ for some $c_1$ and $c_2$ from $K_{-}$.
So one gets $\tilde{T}_3=T_2(u_c : c \in K_{-})$.
In general, the field $\tilde{T}_{n+1}$ is the composite of the fields
$\sigma(T_{n+1})$ with $\sigma \in \Gamma_{n+1}$ and $\sigma(T_{n+1})$
is contained in the field $\sigma(\tilde{T}_n(x_{n+1}))=\tilde{T}_{n}(\sigma(x_{n+1}))$. Again, by applying repeatedly
$\sigma$ to (1) one sees $\wp(\sigma(x_{n+1}))=g(u_{c^{\prime}}+c_n)$ for some
$c^{\prime}\in K_{-}^{n-1}$ and $c_n \in K_{-}$.
\end{proof}
\end{proposition}
Now we shall see that we can restrict to a certain subset of the
$u_c$, namely those for which
$c=(c_1,\ldots,c_{n-2})$ with $c_{n-2} \in \{0,b\}$
for a fixed element $b \neq 0$ of $K_{-}$. For this we note that
$\tilde{T}_3=T_3(u_c)$ for any non-zero $c$ in $K_{-}$. Indeed, given
such $c$ we have the identity
$$
\wp(u_c-x_3+c^2/x_1)= c\, \frac{x_2^{p-1}-1}{x_2^{p-1}+1}
$$
which follows directly from writing out the left hand side. This implies
that
$$
cu_b-bu_c=(c-b)x_3+(bc^2-b^2c)/x_1+ \delta_{c,b}
$$
for some $\delta_{c,b} \in K_{-}$. In general, if for $c^{\prime}=(c_1,\ldots,c_{n-1})$
 and $c=(c_1,\ldots,c_{n-1},\xi)$
we write $u_{c^{\prime}, \xi}$ for $u_c$, we have by a similar argument
 for $c \in K_{-}^{n-1}$
$$
\alpha u_{c^{\prime},\beta}-\beta u_{c^{\prime}, \alpha}= (\alpha -\beta)u_{c^{\prime},0}+
(\beta \alpha^2 -\alpha\beta^2)/(u_{c^{\prime}+ \alpha_{n-1}}) + \eta_{\alpha, \beta}
$$
for some $\eta_{\alpha, \beta} \in K_{-}$ and hence for a fixed $\beta\in K_{-}$
with $\beta \neq 0$ we get
$\tilde{T}_{n+2}=\tilde{T}_{n+1}(u_{c,\xi }: c \in K_{-}^{n-1}, \, \xi \in \{0,\beta
\})$.  We conclude:
\begin{proposition} Let $\beta$ be a non-zero element of $K_{-}$.
The field $\tilde{T}_{n+2}$ is generated over $T_1$ by the set of elements
$\{ u_{c,\xi }: c \in K_{-}^{n-1}, \, \xi \in \{0,\beta\} \}$.
\end{proposition}
In the following we shall also need the following formulas (where we recall
that $\wp(x)=x^p+x$ and $h(x)=(x^{p-1}-1)/(x^{p-1}+1)$).
\begin{lemma}\label{relation}
For all $\alpha, \alpha_1$ in  $K_{-}$ and  $c \in K_{-}^{n-2}$ with $\alpha \neq 0$ we have
\begin{enumerate}
\item $\wp(u_{\alpha}-x_3+\alpha^2/x_1)=\alpha \, h(x_2)$,
\item $\wp(u_{\alpha_1, \alpha}-u_{\alpha , 0}+\alpha^2/(x_2+\alpha_1))=\alpha \, h(u_{\alpha_1})$,
\item $\wp(u_{c,\alpha_{k-1},\alpha}-u_{c,\alpha_{k-1},0}+\alpha^2/(u_{c}+\alpha_{k-1}))=\alpha h(u_{c, \alpha_{k-1}}).$
\end{enumerate}
\end{lemma}
\begin{proof}
This can be proved by direct calculation. As an example we prove the second relation by writing
$$
\begin{array}{rcl}
\wp(u_{\alpha_1, \alpha})&=&g(u_{\alpha_1}+\alpha)=g(u_{\alpha_1})+\alpha\, h(u_{\alpha_1})-\alpha^2\, 1/{\wp(u_{\alpha_1})}\\
\phantom{\wp(u_{\alpha_1, \alpha})}&=&\wp(u_{\alpha_1, 0})+\alpha\, h(u_{\alpha_1})-\alpha^2 \, \wp(1/(x_2+\alpha_1))
\end{array}
$$
and observing that $\wp$ is additive.
\end{proof}
\end{section}
\begin{section}{Splitting Points}
Let $X_n$ (resp.\ $\tilde{X}_n$)
be the irreducible complete smooth algebraic curve defined over $K$
by the function field $T_n$ (resp.\ $\tilde{T}_n$).
Note that $X_1$ is the projective line $\PP^1$.
Here we prove that all the points of the affine line with coordinates not
in $K_{-}$ split completely. At this moment we shall use notations $\pi_n$ and $\tilde{\pi}_n$
for the coverings  $ X_n\rightarrow\PP^{1}$
and $\tilde{X}_n \rightarrow\PP^{1}$ respectively.

\begin{proposition}
Every $K$-rational point of the affine line $\AA^1 \subset \PP^1=X_1$
with coordinate not in $K_{-}$ splits completely
in the tower $\tilde{X}_n$.
\begin{proof}
Since $\tilde{T}_n$ is obtained by adjoining successively the elements
$u_c$ for $c \in K_{-}^{n-2}$ to $\tilde{T}_{n-1}$ we start with a
 $K$-rational point $P=P_1$ not in $K_{-}$
and consider the behavior of points lying over $P$ in these successive
extensions.

Let $P=P_1$ be a point of the affine $x_1$-line with coordinate
$\xi$ in $\FF_{p^2}\backslash K_{-}$. Let ${\rm Nm}$ and ${\rm Tr}$ denote the
trace from $\FF_{p^2}$ to $\FF_{p}$. By the identity
$$
\frac{\xi^{p+1}}{\xi^p+\xi}= \frac{{\rm Nm}(\xi)}{{\rm Tr}(\xi)}
\eqno(2)
$$
it is clear that this expression lies in $\FF_p^*$ and
it is immediate that $P$ splits completely in the field extension $T_2/T_1$
given by adjoining a root of $Y^p+Y={\rm Nm}(\xi)/{\rm Tr}(\xi)$
and the $x_2$-coordinate of any point $P_2$ over $P$
has coordinate $\eta$ in ${\FF}_{p^2}\backslash K_{-}$. So we can repeat the argument and see that
$P_2$ splits completely in the extension $T_3/T_2$. Since for $c \in K_{-}$ $u_c$
is a root of $X^p+X={\rm Nm}(\eta+a)/{\rm Tr}(\eta)$ and the right hand side lies
in $\FF_{p}^*$ we see again that $P$ splits completely in $\tilde{T}_3$.
For the general step we need the following lemma.

\begin{lemma}
Let $Q$ be a point on $\tilde{X}_n$ lying over
$P\in \AA^1(\FF_{p^2})\backslash
\AA^1(K_{-})$. Then for $n \geq 3$, any $c \in K_{-}^{n-2}$ and $\alpha \in K_{-}$
the value $u_c(Q)$ lies in
$\FF_{p^2}\backslash K_{-}$  the polynomial
$$
X^p+X-\frac{(u_c+\alpha)^{p+1}}{u_c^p+u_c}
$$
splits completely into linear factors over $\FF_{p^2}$ at the point $Q$.
\begin{proof}
We use induction on $n$ starting with $n=3$. For $n=3$ relation (2) shows
that $g(x_2+\alpha)(Q)$ lies in $\FF_p^*$ and hence $u_a(Q) \in \FF_{p^2}\backslash K_{-}$
for any $\alpha \in K_{-}$. We denote by $Q_{n-1}$ the image of $Q$ on $\tilde{X}_{n-1}$.
Assume now that $u_c(Q_{n-1}) \in \FF_{p^2}\backslash K_{-}$ for an arbitrary $c \in K_{-}^{n-3}$.
So for any $\alpha \in K_{-}$ the expression ${\rm Nm}(u_c+\alpha)/{\rm Tr}(u_c^p+u_c)$
is $\FF_p$-valued and does not vanish in $Q_{n-1}$.
But that implies that our polynomial $f_{c, \alpha}$ evaluated at
$Q_{n-1}$ factors linearly and has no roots in $K_{-}$. It follows that
$u_{c, \alpha}(Q) \in \FF_{p^2}\backslash K_{-}$.
\end{proof}
\end{lemma}
\end{proof}
\end{proposition}

\begin{corollary}
The finite field $\FF_{p^2}$ is the full constant field of the function field of the curve $\tilde{X}_n$.
\end{corollary}

\begin{corollary}
The curve $\tilde{X}_n$ has at least $(p^2-p) [\tilde{T}_n : T_1]$ $\FF_{p^2}$-rational
points.
\end{corollary}
\end{section}
\begin{section}{Ramification over Zero}
In this section we calculate the contribution to the different
of the ramifying points of $\tilde{X}_n$ which lie over the point $P_0$
of $X_1$ given by $x_1=0$.

\begin{proposition}\label{path} Let $n \geq 4$. There exists  points $Q_i
\in \tilde{X}_i$ for $i=1,\ldots, n$ such that $Q_{i+1}$ lies over
$Q_i$ and such that $Q_3$ is unramified over $Q_1$ and such that for $i\geq 3$
the point $Q_{i+1}$ ramifies over $Q_i$ with ramification index $e=p$ and different degree $d=2(p-1)$.
\begin{proof}
For the proof we first observe that if $Q_3^{\prime}\in X_3(\FF_{p^2})$
is the point of $X_3$
defined by $x_1=x_2=x_3=0$ then $Q_3^{\prime}|Q_1$ is unramified as follows from
\cite{G-S2}.
Moreover, $\tilde{T}_3$ can be generated over $T_3$ by adjoining a root
of $T^p+T=c \, h(x_2)$ for an arbitrary $c \in K_{-}^*$.
The right-hand side of this has value $-c$ in $Q_3^{\prime}$,
so $Q_3^{\prime}$ is inert giving a point $Q_3$ on $\tilde{X}_3$. We also
observe for later use that for any $c \in K_{-}$ the function $u_c+c^2/x_1$
is regular at $Q_3$.

The proof of the proposition is now by induction starting
with the case $n=3$ just settled. We assume having established the existence of a point $Q_{n+2}$
on $\tilde{X}_{n+2}$ satisfying the following properties $P(n+2)$.
We shall denote the zero vector with i coordinates by $0^i$.
\begin{property}
We say that a point $P$ of the curve $\tilde{X}_{n+2}$ has property $P(n+2)$ if the following conditions hold
\begin{enumerate}
\item The point $P$ is a zero of  the functions $x_1, x_2, \ldots, x_n$.
\item For any $\alpha \in K_{-}^{*}$ and non-zero $c \in K_{-}^{n}$
function $u_{c}$ has a pole at the point $P$ and this pole is simple if $c=(\alpha, 0^{n-1})$.
\item For any $\alpha, \beta \in K_{-}^{*}$ and $c \in K_{-}^{n-1}$,
the function  $u_{\beta, c}-\beta^2/\alpha^2 u_{\alpha, 0^{n-1}}$ is regular at $P$.
\item For any $m \geq 2$ and any element $c=(c_1, \ldots, c_m) \in K_{-}^{*}\times K_{-}^{m-1}$ the function
$u_{0^{n-m}, c}-u_{0^{n-m-1}, c_1,  0^{m-2}}$ is regular at $P$.
\item The function $u_{c}+ c_n^2/x_n$ is regular at $P$ where $c=(0, \ldots, 0, c_n) \in K_{-}^{n}$.
\end{enumerate}
\end{property}
We construct  the field extension $\tilde{T}_{n+3}$ of $\tilde{T}_{n+2}$ by
successively adjoining elements $u_c$ with $c=(c_1, \ldots, c_{n+1}) \in K_{-}^{n+1}$ and we analyze when
we get contribution to the different from ramification. For  a convenience analysis
we separate our indices into the following sorts:
\begin{enumerate}
\item $c_1=\ldots =c_{n+1}=0$,
\item $c_1=\ldots =c_{n}=0$, $c_{n+1}\neq 0$,
\item $c_1\neq 0$, $c_2=\ldots =c_{n+1}=0$,
\item $c_1\neq 0$, $c_{n+1}=0$,
\item $c_1\neq 0$, $c_{n+1}\neq 0$,
\item $c_1=\ldots =c_s=0$, $c_{s+1}\neq 0$, $c_{n+1}=0$, $n+1-s\geq 2$ and $s\geq1$,
\item $c_1=\ldots =c_s=0$, $c_{s+1}\neq 0$, $c_{n+1}\neq 0$, $n+1-s\geq 2$ and $s\geq1$,
\item $c_1=0 \ldots=c_{n-1}=0$, $c_{n} \neq 0$ and $c_{n+1}=0$,
\item $c_1=0 \ldots=c_{n-1}=0$, $c_{n} \neq 0$ and $c_{n+1}\neq 0$;
\end{enumerate}
and adjoin successively elements $u_c$ with $c=(c_1,\ldots,c_{n+1})$ of these types.

We shall show that only elements of type 3) contribute to the different.
We select the generators and their polynomials at each stage such that we are able to apply
Artin-Schreier reduction(\cite{N-X} Proposition 3.1.10 p.64) or Kummer's theorem(\cite{S} Theorem III.3.7 p.76).
As is well-known a polynomial of the form $T^p+T+w \in F[T]$, with $F$ a field extension of $\FF_{p^2}$, is
either irreducible or splits into linear factors. If such a polynomial is reducible, then adjoining a root
will not extend the field, but this will not lead to confusion.
\par
    If we adjoin an element of type 1) then we get a function field $F_1=\tilde{T}_{n+2}(u_c)$
with $c=(0,\ldots,0)$ and this extension is actually generated  by an element $x_{n+3}$
satisfying an equation $\wp(x_{n+3})=g(x_{n+2})$. Now observe that
$g(x_{n+2})$ vanishes at $Q_{n+2}$, so $Q_{n+2}$ splits in this extension
giving us a point $Q_{n+2,1}$ on the corresponding curve, such that $x_{n+3}$ vanishes at $Q_{n+3}$.
(In the next section we shall show that the polynomial $\wp(X)-g(x_{n+2}) \in \tilde{T}_{n+2}[X]$ is in fact irreducible.)

     Next we treat the case (2), where we adjoin a root of $u_c$ of $\wp(u_c-x_{n+3}+c_{n+1}^2/x_{n+1})=c_{n+1} \, h(x_{n+2})$
for $c=(0,\ldots,0,c_{n+1})$ with $c_{n+1}\neq 0$. Now $h(x_{n+2})$ assumes the
value $-1$ at $Q_{n+2,1}$. So in this extension $F_1(u_c)$ the point
$Q_{n+2,1}$ does not ramify, giving us the point $Q_{n+2,1}^{\prime}$ such that $u_c+c_{n+1}^2/x_{n+1}$ is
a regular function at the point $Q_{n+2,1}^{\prime}$.

If we adjoin repeatedly such elements $u_c$ with $c$ of type 2) 
then Abhyankar's
lemma(cf., \cite{S} Prop.\ III.8.9) implies that in the composite 
field $F_2:=F_1(\{u_c: \hbox{\rm c of type 2)}\})$
there is a point $Q_{n+2,2}$ lying over $Q_{n+2,1}$ with $e(Q_{n+2,2}|Q_{n+2,1})=1$ and
the functions $u_c+c_{n+1}^2/x_{n+1}$ are regular at the point $Q_{n+2,2}$ for any $c$ of type 2).

    For the rest of the proof we fix an element $\alpha \in K_{-}^{*}$.
Now let $F_3$ be the field obtained by adjoining to $F_2$ one element $u_c$ with $c_1=\alpha$
and $c_2=\ldots =c_{n+1}=0$. Then the defining equation for this field extension is $f_c=0$. Since we know
that $g(u_c)$ has a simple pole at $Q_{n+2,2}$ the point $Q_{n+2,2}$ ramifies
giving one point $Q_{n+2,3}$ with contribution $2(p-1)$ to the
different and the function $u_{\alpha, 0^{n}}$ has a simple pole at point $Q_{n+2,3}$.

    If we adjoin an element of type 4) then combining the two relations $\wp(u_c)=u_{c^{\prime}}+\cO(1/u_{c^{\prime}})$ and
$\wp(u_{\alpha, 0^n})=u_{\alpha, 0^{n-1}}+\cO(1/u_{\alpha,
0^{n-1}})$ at the point $Q_{n+2,3}$, we see that the point
$Q_{n+2,3}$ does not ramify and gives a point $Q_{n+2,3}^{\prime}$
such that the function $u_c-c_1^2/\alpha^2 u_{\alpha, 0^n}$ is
regular at  the point $Q_{n+2,3}^{\prime}$. (Here and in the following the
symbol $\cO(t)$ at a point $P$ means a function of the form $ut$ with $u$ is a regular
function at $P$).

   If we adjoin repeatedly such elements $u_c$ with $c$ of type 4), then Abhyankar's lemma implies that
in the composite field $F_4:=F_3(\{u_c: \hbox{\rm c of type 4)}\})$ there is a point $Q_{n+2,4}$ lying
over $Q_{n+2,3}$ with $e(Q_{n+2,4}|Q_{n+2,3})=1$ and such that  $u_c-c_1^2/\alpha^2 u_{\alpha, 0^n}$
is the regular function at the point $Q_{n+2,4}$ for any element $c$ of type 4).

   Next we are going to adjoin the elements of type 5). Observe that $h(u_{c^{\prime}})$ has value $1$
at the point $Q_{n+2,4}$. Therefore in view of the third relation in
Lemma \ref{relation}, we obtain that the point $Q_{n+2,4}$ does not ramify, giving us a point $Q_{n+2,4}^{\prime}$
on the corresponding curve such that the function $u_c-c_1^2/\alpha^2 u_{\alpha, 0^n}$ is regular
at the point $Q_{n+2,4}^{\prime}$.

        If we adjoin repeatedly such elements $u_c$ with $c$ of type 5) then Abhyankar's
lemma implies that in the composite field $F_5:=F_4(\{u_c: \hbox{\rm c of type 5)}\})$
there is a point $Q_{n+2,5}$ lying over $Q_{n+2,4}$ with $e(Q_{n+2,5}|Q_{n+2,4})=1$ and such that the functions
$u_c- c_1^2/\alpha^2 u_{\alpha 0^n}$ are regular at the point $Q_{n+2,5}$.

If we adjoin an element of type 6), then we write 
$\wp(u_c)=u_{c^{\prime}}+\cO(1/u_{c^{\prime}})$ and 
$\wp(u_{0^{s},c_{s+1},0^{n-s}})=
u_{0^{s-1},c_{s+1},0^{n-s-1}}+\cO(1/u_{0^{s-1},c_{s+1},0^{n-s-1}})$
at the point $Q_{n+2,5}$, to get  a point $Q_{n+2,5}^{\prime}$ that lies over $Q_{n+2,5}$ with
$e(Q_{n+2,5}^{\prime}|Q_{n+2,5})=1$ and such that the function $u_c- u_{0^{s},c_{s+1},0^{n-s}}$ is regular
at the point $Q_{n+2,5}^{\prime}$.

In case we adjoin repeatedly such elements $u_c$ with $c$ 
of type 6), Abhyankar's
lemma implies that in the composite field 
$F_6:=F_5(\{u_c: \hbox{\rm c of type 6)}\})$
there is a point $Q_{n+2,6}$ lying over $Q_{n+2,5}$ with $e(Q_{n+2,6}|Q_{n+2,5})=1$ and
such that the functions $u_c-u_{0^{s},c_{s+1},0^{n-s}}$ are regular at the point $Q_{n+2,6}$.

     If we adjoin an element $u_c$ with $c$ of type 7), then in view of relation 3) of Lemma \ref{relation}
and such that the value of $h(u_{c^{\prime}})$ equals $1$ at the point  $Q_{n+2,6}$ we get a point $Q_{n+2,6}^{\prime}$
lying over $Q_{n+2,6}$ such that $e(Q_{n+2,6}^{\prime}|Q_{n+2,6})=1$ and  the function
$u_c-u_{c^{\prime},0}$  is  regular at the point $Q_{n+2,6}^{\prime}$, so $u_c- u_{0^{s},c_{s+1},0^{n-s}}$
is also a regular function at the point $Q_{n+2,6}^{\prime}$.
Now Abhyankar's lemma applied to the composite field $F_7=F_6(\{u_c: \hbox{\rm of type 7)}\})$
shows that there exists a point $Q_{n+2,7}$ lying over
$Q_{n+2,6}$ such that $e(Q_{n+2,7}|Q_{n+2,6})=1$ and such that the function $u_c-u_{c^{\prime},0}$
 is  regular at the point $Q_{n+2,7}$, hence $u_c- u_{0^{s},c_{s+1},0^{n-s}}$ is  regular at the point $Q_{n+2,7}$.

     The last two steps differ little from the previous ones.
First we adjoin an element $u_c$ with $c$ of type 8). The relation
$\wp(u_c+c_{n}^2/x_{n-1})=u_{c^{\prime}}+c_{n}^2/x_{n}+\cO(1/{u_c^{\prime}})+\cO(x_n)$ at the point $Q_{n+2,7}$,
implies that $Q_{n+2,7}$ does not ramify, giving us a point $Q_{n+2,7}^{\prime}$ on the corresponding curve
such that the function $u_c+c_{n}^2/x_{n-1}$  is regular at the point $Q_{n+2,7}^{\prime}$
and hence $u_c-u_{0^{n-2},c_{n}}$is also regular.
Applying Abhyankar's lemma to the composite field $F_8=F_7(\{u_c: \hbox{\rm c of type 8)}\})$  we see that
there is a point $Q_{n+2,8}$ lying over $Q_{n+2,7}$ with $e(Q_{n+2,8}|Q_{n+2,7})=1$
and such that the functions  $u_c-u_{0^{n-1},c_{n}}$ are  regular  at the point $Q_{n+2,8}$.

     Finally we adjoin the elements $u_c$, where $c$ is type 9). Since the function
$\wp(u_c-u_{0^{n-1}, c_{n}, 0}+c_{n+1}^2/(x_{n+1}+c_{n})$ is  regular at the point $Q_{n+2,8}$
(the value $h(u_{0^{n},c_{n}})$ is $1$ at this point) we see that the point $Q_{n+2,8}^{\prime}$ does not ramify
with   $u_c-u_{0^{n-1}, c_n, 0}$  being a regular function
at the point $Q_{n+2,8}^{\prime}$ and hence $u_c-u_{0^{n-2}, c_n}$ is  regular as well.
Abhyankar's lemma applied to the composite field $F_9=F_8(\{u_c: \hbox{\rm c of type 9)}\})$ shows
the existence of a point $Q_{n+2,9}$( which is now called $Q_{n+3}$), such that $e(Q_{n+2,9}|Q_{n+2,8})=1$
and the functions $u_c-u_{0^{n-2}, c_n}$ are regular at the point $Q_{n+2,9}=Q_{n+3}$.

    So we conclude that $F_9=\tilde{T}_{n+3}$ and our proof is finished  but for the remark
that Property (2) of $P(n+3)$ holds because  the function $u_c$ has a pole at $Q_{n+3}$ for non-zero $c$ by the induction hypothesis and
since function $\wp(u_{c^{\prime},c_{n+1}})$ has a pole at $Q_{n+3}$.

\end{proof}
\end{proposition}

    Since $\tilde{X}_n$ is a Galois covering of $X_1$, for calculating the contribution to the ramification divisor
of all points lying over $P_0$ it suffices to calculate the contribution of one such point.
This contribution was calculated in Proposition \ref{path}. Collecting results we obtain the following corollary.

\begin{corollary}
For $n\geq 4$ let $D_{n}$ be the divisor on the curve  $\tilde{X}_{n}$ such that
$v_P(D_n)=v_P({\rm Diff}(\tilde{T}_n / T_1))$ and $P \cap T_1=P_0$   for any $P \in {\rm Supp}(D_n)$.
Then we have
$$
\deg(D_n)=2\deg(\tilde{\pi}_n)(1-p^{3-n}).
$$
\begin{proof}
For any point $Q$ on curve $\tilde{X}_{n}$ lying over $P_0$, we have $d(Q|P_0)=d(Q_{n}|P_0)$
with $Q_n$ the point mentioned in Proposition \ref{path}. We obtain
$$
\begin{array}{rl}
\deg(D_n)&=(1+ \ldots+ p^{n-4})2(p-1)\# \{Q \in \tilde{X}_n: Q|P_0\}\deg(Q_n)\\
\phantom{\deg(D_n)}&=2\deg(\tilde{\pi}_n)(1-p^{3-n}).
\end{array}
$$
\end{proof}
\end{corollary}
\end{section}
\begin{section}{The other rational points}
In this section we calculate the contribution to the different of the ramifying points of $\tilde{X}_n$
which lie over the point $P_1$ equal to $\infty$ or $a \in K_{-}^{*}$ in $\PP^1$ .
These two kinds of points have the same behavior.
\begin{proposition}\label{path2} Let $P_1$ be a rational point on $\PP^{1}(K)$ with
coordinate   $a \in K_{-}^{*}\cup \infty$. Then there exists points $P_i$ on $\tilde{X}_i$ for $i=1,\ldots, n$ such that $P_{i+1}$ lies over
$P_i$ and such that the point $P_{i+1}$ ramifies over $P_i$
with ramification index $e=p$ and different degree $d=2(p-1)$.
\begin{proof}
In view of the fact that the function $g(x_1)$ has a simple pole at $P_1$ we get that $P_1$ ramifies, giving us
a point $P_2$ on curve $X_2$ with $d(P_2|P_1)=2(p-1)$ and
such that the function $x_2$ has a simple pole at $P_2$. Therefore  the function $g(x_2)$ again has a simple pole
at the point $P_2$ and the point $P_2$ ramifies, yielding a point $Q$ on curve $X_3$
with $d(Q|P_2)=2(p-1)$ and such that the function $x_3$ has a simple pole at $Q$.
Now to reach the curve $\tilde{X}_3$ we shall adjoin an element $u_{c}$ with $c \in K_{-}^{*}$;
in this case the function $h(x_2)$ has the value $1$ at the point $Q$, so the point $Q$ does not ramify, leading  us
to a point $P_{3}$ on curve $\tilde{X}_3$ and such that the function $u_c-x_3$ is regular,
hence the function $u_c$ has a simple pole at the point $P_3$.

The proof of the proposition is analogous to that of Proposition \ref{path}, starting with case $n\leq 3$ just settled.
We distinguish various cases depending on $c \in K_{-}^{n+1}$ and the property $P(n)$ is replaced by the Property $S(n)$ below.
We construct the field extension
$\tilde{T}_{n+3}$ over $\tilde{T}_{n+2}$ by successively adjoining elements $u_{c}$ with $c \in K_{-}^{n+1}$.
Like in the previous  section we distinguish several cases:
\begin{enumerate}
\item $c_1=\ldots =c_{n+1}=0$,
\item $c=(c^{\prime},0)$, with $c^{\prime} \in K_{-}^{n}$,
\item $c=(c^{\prime},c_{n+1})$, with $c^{\prime} \in K_{-}^n$ and $c_{n+1}\in K_{-}^{*}$.
\end{enumerate}
\begin{property}\label{property_S}
We say that a point $P$ of the curve $\tilde{X}_{n+2}$ has property $S(n+2)$ if the following conditions hold
\begin{enumerate}
\item The function $u_c$ has a simple pole at the point $P$ for any $c \in K_{-}^{n}$.
\item The function $u_c-u_{0^{n}}$ is regular at the point $P$ for any $c \in K_{-}^{n}$.
\end{enumerate}
\end{property}
If we adjoin an element of type (1) then we obtain a function field $F_1:=\tilde{T}_{n+2}(u_c)$. Now observe
that the function has a simple pole at the point $P_n$, therefore the polynomial $\wp(X)-g(x_{n+2}) \in \tilde{T}_{n+2}[X]$
is irreducible and $P_{n+2}$ ramifies, and provides us with a point
$P_{n+2,1}$ on the corresponding curve such that the function $u_c$ has a simple pole at $P_{n+2,1}$.
(In particular, we have proved irreducibility of the polynomial $\wp(X)-g(x_{n+2})$ over $ \tilde{T}_{n+2}$
 as we promised  in the previous section.)

Next we treat the case 2) and adjoin $u_c$, which is  a root of $f_c=0$. If we combine the two relations
$\wp(u_c)=u_{c^{\prime}}+\cO(1/u_{c^{\prime}})$ and $\wp(u_{0^{n+1}})=u_{0^{n}}+\cO(1/u_{0^{n}})$ at the point $P_{n+2,1}$
then we obtain that $P_{n+2,1}$ does not ramify, yielding a point $P_{n+2,1}^{\prime}$ such that
the function $u_c-u_{0^{n+1}}$ is regular at the point $P_{n+2,1}^{\prime}$.
If we adjoin repeatedly such elements $u_c$ with $c$ of type 2) then Abhyankar's
lemma implies that in the composite field $F_2:=F_1(\{u_c: \hbox{\rm c of type 2)}\})$
there is a point $P_{n+2,2}$ lying over $P_{n+2,1}$ with $e(P_{n+2,2}|P_{n+2,1})=1$ and such that the functions
$u_c-u_{0^{n+1}}$ are regular at the point $P_{n+2,2}$.

If we adjoin an element $u_c$ with $c$ of type 3), then in view of relation 3) of Lemma \ref{relation}
and the fact that the value of $h(u_{c^{\prime}})$ equals $1$ at the point  $P_{n+2,2}$ we get a point $P_{n+2,2}^{\prime}$
lying over $P_{n+2,2}$ such that $e(P_{n+2,2}^{\prime}|P_{n+2,2})=1$ and such that the function
$u_c-u_{c^{\prime},0}$  is  regular at the point $P_{n+2,2}^{\prime}$;
therefore $u_c-u_{0^{n+1}}$ is also a regular function at the point $P_{n+2,2}^{\prime}$.
Finally, Abhyankar's lemma applied to the composite field $F_3=F_2(\{u_c: \hbox{\rm of type 3)}\})$
shows that there exists a point $P_{n+2,3}$(which we call $P_{n+3}$) lying over
$P_{n+2,2}^{\prime}$ such that $e(P_{n+2,3}|P_{n+2,2}^{\prime})=1$ and  the functions $u_c-u_{c^{\prime},0}$
are  regular at the point $P_{n+2,2}$; hence $u_c- u_{0^{n+1}}$ is regular at the point $P_{n+2,3}$ for all such $c$.
To finish our proof, we remark that the first condition of Property(\ref{property_S}) follows from the strict triangle
inequality for discrete valuations applied to the relation $\wp(u_c)=g(u_{c^{\prime}}+c_{n+1})$.
\end{proof}
\end{proposition}
Since $\tilde{X}_n$ is a Galois covering of $X_1$, for calculating the contribution to the ramification divisor
of all points lying over $P_1$ it suffices to calculate the contribution of one such point.
This contribution was calculated in Proposition \ref{path2}. We thus obtain the following corollary.
\begin{corollary}
For $n\geq 5$ let $L_{n}$ be the divisor on the curve  $\tilde{X}_{n}$ such that
$v_P(L_n)=v_P({\rm Diff}(\tilde{T}_n / T_1))$ and $P \cap T_1$ is the point $P_1$ for any $P \in {\rm Supp}(L_n)$.
Then we have
$$
\deg(D_n)=2(p-p^{2-n})\deg(\tilde{\pi}_n).
$$
\begin{proof}
For any point $P$ on curve $\tilde{X}_{n}$ lying over $P_1$, then $d(P|P_1)=d(P_{n}|P_0)$.
For the point $P_n$ from Proposition \ref{path2} we have
$d(P_n|P_1)=(1+ \ldots+ p^{n-4})2(p-1)+2p^{n-3}(p^2-1)$. As a result we obtain
$$
\deg(D_n)=\frac{\deg(\tilde{\pi}_n)}{e(P_1)f(P_1)}2(p^{n-1}-1)\deg(P_n)=\frac{2(p^{n-1}-1)}{p^{n-2}}\deg(\tilde{\pi}_n).
$$
\end{proof}
\end{corollary}
\end{section}
\begin{section}{The genus of the curves and the optimality of 
the new sequence}

In this section we shall show that our sequence of curves attains
the Drinfeld-Vl\u{a}du\c{t} bound. We show that $\lim_{n
\rightarrow \infty}N(\tilde{X}_n)/g(\tilde{X}_n)=p-1$. Since we
already estimated  the number of rational point of the new curves,
we only need to calculate the genus of these curves. We are going
to show that the different of the covering $\tilde{X}_n$ over
$\PP_1$ is the sum of two divisors $D_n$ and $L_n$ described in
section 4 and 5 respectively. After that the calculation of the
genus is simple.
\begin{proposition}
Let $n\ge 5$ we have that ${\rm Diff}(\tilde{X}_n/\PP^{1})=D_n+L_n$ with the divisors $D_n$ and $L_n$
defined in sections 4 and 5, respectively.
\begin{proof}
By Artin-Schreier reduction we see that if a point of the curve $\tilde{X}_{n+2}$ contributes to the different
of the covering $\tilde{X}_{n+3}/\tilde{X}_{n+2}$ then it is a pole of the function $g(u_{a}+a_{n+1})$
with $(a,a_{n+1}) \in K_{-}^{n+1}$. To find all these points we shall consider the divisors of the functions $u_{a}+a_{n}$
with $(a,a_{n+1}) \in K_{-}^{n+1}$.
We need a lemma.
\begin{lemma}\label{supp}
Let a point $P$ of the curve $\tilde{X}_{n+2}$ be either a pole or a zero of some function $u_{a}+a_{n}$
with $(a,a_{n+1}) \in K_{-}^{n+1}$. Then it lies over rational point on $\PP^{1}$ with coordinate in
$K_{-}\cup \{\infty\}$.
\begin{proof}
We proceed by induction on $n$.
Let $n=1$, i.e. we have a point $P$ on the curve $\tilde{X}_{3}$ which is either a zero or a pole of the function $u_{a_1}+a_{2}$.
In this case it has to be either a zero or a pole of the function $g(x_2+a_1)$ and hence $P\in {\rm supp}(x_2+\gamma)$ for
some $\gamma \in K_{-} $. Since $\wp(x_2+\gamma)=g(x_1)$ we get that $P$ is either a zero or a pole of the function $x_1-\beta$
for some $\beta \in K_{-}$.

For the general case let $P$ be a point of the curve $\tilde{X}_{n+2}$ in ${\rm supp}(u_a+a_{n+1})$
with  $(a,a_{n+1}) \in K_{-}^{n+1}$.
Then from the relation $\wp(u_a+a_{n+1})=g(u_{a^{\prime}+a_{n}})$
we get that $P$ in ${\rm supp}(u_{a^{\prime}+a_{n}})$ as well. Using the induction hypothesis we obtain
that $P$ lies over some rational point on $\PP^{1}$ with coordinate in $K_{-}\cup \{\infty\}$.
\end{proof}
\end{lemma}
Now we show that if a point $P$ of the curve $\tilde{X}_{n+2}$ contributes
to the different ${\rm Diff}(\tilde{X}_{n+3}/\tilde{X}_{n+2})$ then it is a pole
of the function $g(u_a+a_{n+1})$ with  $(a,a_{n+1}) \in K_{-}^{n+1}$ and hence
it is either a pole or a zero of the function $u_a+a_{n+1}$. Therefore by Lemma \ref{supp} the point $P$
lies over rational point on $\PP^{1}$ with coordinate in $K_{-}\cup \{\infty\}$.
Applying that
$
{\rm Diff}(\tilde{X}_{n+3}/\PP^{1})={\rm Diff}(\tilde{X}_{n+3}/\tilde{X}_{n+2})+
{\rm (\pi_{n+3,n+2})_{*}}({\rm Diff}(\tilde{X}_{n+2}/\PP^{1}))
$
with $ {\rm \pi_{n+3,n+2}}:\tilde{X}_{n+3}\rightarrow\tilde{X}_{n+2}$
and  the induction hypothesis we finish our proof.
\end{proof}
\end{proposition}
Next we calculate the genus of the curve $\tilde{X}_n$.
\begin{corollary}
For $n >4$ the genus of the curve $\tilde{X}_{n}$
is given by the formula\\
\centerline{$g(\tilde{X}_n)=\deg(\tilde{\pi}_n)(p-p^{3-n}-p^{2-n})+1$.}
\begin{proof}
By the Hurwitz genus formula for the covering $\tilde{\pi}_n$ we have
$$
\begin{array}{rl}
g(\tilde{T}_n)&=1/2\deg({\rm Diff}(\tilde{\pi}_n))-\deg(\tilde{\pi}_n)+1\\
&=1/2(\deg(D_n)+\deg(L_n))-\deg(\tilde{\pi}_n)+1\\
&=\deg(\tilde{\pi}_n)(p-p^{3-n}-p^{2-n})+1.
\end{array}
$$
\end{proof}
\end{corollary}
Since we know the genus of the curve $\tilde{X}_n$  we can now present the main result of this article, namely
 that the new sequence of curves is optimal.
\begin{theorem}
The sequence of curves $\{\tilde{X}_{n}\}_{n\ge 1}$  attains 
the Drinfeld-Vl\u{a}du\c{t} upper-bound,
i.e., 
$\lim_{n \rightarrow \infty}N(\tilde{X}_n)/g(\tilde{X}_n)=p-1$.
\begin{proof}
We have
$$
\lambda(\{\tilde{X}_{n}\})=\lim_{n \rightarrow \infty}\frac{N(\tilde{X}_n)}{g(\tilde{X}_n)} \geq
\lim_{n \rightarrow \infty}\frac{\deg(\tilde{\pi}_n)(p^2-p)}{\deg(\tilde{\pi}_n)(p-p^{3-n}-p^{2-n})+1}=p-1.
$$
Therefore the Drinfeld-Vl\u{a}du\c{t} upper bound provides the equality.
\end{proof}
\end{theorem}

At the end we would like to remark that the result can be generalized
with $p$ replaced by any power of an odd prime. For this one should 
use Kummer's theorem instead of
Artin-Schreier reduction
for proving that points do not contribute to the different.
It is based on the fact that the polynomial $f_c(X)$, 
with $p$ changed to $q$, gives us a separable
polynomial under reduction at certain points,
hence its irreducible factors also give us separable polynomials 
under reduction at such points.
Therefore those points are unramified by Kummer's theorem.
\end{section}

\end{document}